\RequirePackage{ifpdf}
\ifpdf 
\documentclass[pdftex]{sigma}
\else
\documentclass{sigma}
\fi

\begin{document}

\allowdisplaybreaks

\renewcommand{\thefootnote}{$\star$}

\renewcommand{\PaperNumber}{063}

\FirstPageHeading

\ShortArticleName{Exterior Dif\/ferential Systems for Yang--Mills Theories}

\ArticleName{Exterior Dif\/ferential Systems\\ for Yang--Mills Theories\footnote{This paper is a
contribution to the Special Issue ``\'Elie Cartan and Dif\/ferential Geometry''. The
full collection is available at
\href{http://www.emis.de/journals/SIGMA/Cartan.html}{http://www.emis.de/journals/SIGMA/Cartan.html}}}

\Author{Frank B. ESTABROOK}

\AuthorNameForHeading{F.B. Estabrook}

\Address{Jet Propulsion Laboratory, California Institute of Technology, Pasadena, CA 91109, USA}
\Email{\href{mailto:frank.b.estabrook@jpl.nasa.gov}{frank.b.estabrook@jpl.nasa.gov}}

\ArticleDates{Received July 18, 2008, in f\/inal form September 03,
2008; Published online September 12, 2008}

\Abstract{Exterior dif\/ferential systems are given, and their Cartan characters calculated,
for Maxwell and $SU(2)$-Yang--Mills equations in dimensions from three to six.}


\Classification{58A15}

\bigskip

\rightline{\it In memory of my friend and colleague Hugo D.~Wahlquist [1930--2008]}

\section{Introduction}
In this work Exterior Dif\/ferential Systems (EDS) are formulated for abelian and $SU(2)$-gauge theories in 3 to 6 dimensions, beginning with Maxwell vacuum electrodynamics and SU(2)-Yang--Mills theory in 4 dimensions.  The generalization to larger groups and $n$ dimensions becomes apparent.  The results are new in that the Cartan integer character tables are explicitly calculated, showing the EDS' to be well posed f\/ield theories. Cartan characters yield the numbers of evolution equations, of f\/irst order (second order, etc.), of constraint equations (or integrability conditions), and of gauge degrees of freedom;  well posed EDS' allow systematic derivation of isovectors and similarity solutions, conservation laws and prolongations.

An EDS for vacuum Maxwell theory in 4 dimensions has been set on a 14 dimensional manifold \cite{EHW, Her, FECons} with variables $A_{i}$, $F_{ij} = -F_{ji}$, $x^{i}$, $i=1,2,3,4 $.  Its generalization to Yang--Mills $SU(2)$, in terms of 34 variables $A^{a}_{i}$, $F^{a}_{ij} = -F^{a}_{ji}$, $x^{i}$, $i=1,2,3,4$, $a=1,2,3 $ was also partially formulated in the f\/irst two of these references. An analogous EDS for vacuum general relativity was more recently found, and Cartan characters and conservation laws calculated for both it and vacuum Maxwell theory~\cite{FEMath,FECons}.  These EDS' are generated by sets of ``torsion'' 2-forms and dynamic 3-forms, and these factor the terms of a Cartan--Poincar\'{e} (or multisymplectic) 5-form $d\Lambda$ expressing the underlying variational principle.  The  generalized EDS' in n dimensions to be given in Sections~\ref{sec2} and~\ref{sec3} will follow this structure, using torsion 2-forms and $(n-1)$-forms.

It should be noted that this approach to f\/inding a Cartan $n$-form $\Lambda$, or Lagrangian density, for a f\/ield theory dif\/fers from one that constructs a Cartan--Poincar\'{e} form $d\Lambda$ using ``contact'' 1-forms and dynamic $n$-forms~\cite{Hoj,Bet}.  That approach has been carefully worked out as ``multicontact'' geometry~\cite{Bry}.  In 4-dimensional Maxwell theory, e.g., the multicontact EDS is set on a 24 dimensional manifold, spanned by four potential f\/ields $A^{i}$, 16 f\/ields $A^{i}_{.j}$ and the $x^{i}$.  The Cartan characteristic integer table is quite dif\/ferent.

The Cartan forms and EDS' in this work were manipulated and checked algebraically \mbox{using} the MATHEMATICA suite of programs AVF (Algebra Valued Forms) written by the late H.D.~Wahlquist~\cite{Wahl}.  The integer characters were calculated using his Monte Carlo program \mbox{characters.m}, which assigns random integer values to unknown components of successively calculated vectors spanning integral manifolds of an EDS, vectors satisfying the nested linear equations of Cartan--Kuranishi integration theory.  After three or so of these are calculated, using integer arithmetic, the new vector components solved for become rational numbers with scores of di\-gits, and the consequent calculation of the ranks of the resulting matrices (and so the integer characters) is absolutely robust. (In principle of course this program could be rewritten to use the f\/irst few hundred primes instead of random integers, making the calculation of characters exact.)

In 2003 Wahlquist circulated on CD a f\/inal version of his principal MATHEMATICA notebooks, with documentation and demonstrations, as ``AVF Programs''.
This is now available from the SIGMA web-site: 
\url{http://www.emis.de/journals/SIGMA/2008/063/Wahlquist2003.zip}.

\section[Maxwell and Yang-Mills]{Maxwell and Yang--Mills}\label{sec2}

We will be generalizing the EDS and Cartan form for the vacuum Maxwell equations;  these are most conveniently written using a 2-form $F=\tfrac{1}{2} F_{ij} dx^{i} \wedge dx^{j}$ and a dual 2-form $*F= \tfrac{1}{4} F_{ij} dx^{k} \wedge dx^{l} \epsilon^{ij}_{..kl}$ involving the f\/ields and a 1-form $A=A_{i}dx^{i}$ coding the potentials.  $i,j = 1,\dots,4$, and indices are raised/lowered with the Lorentz metric.  The EDS is generated~\cite{FECons} by a~2-form $\theta = dA-F$, its closure $d\theta$ and another 3-form $\psi = -d*F$.  The Cartan characters are $s_{0}=0$, $s_{1}=1$, $s_{2}=3$, $s_{3}=5$, $s_{4}=1$.  The genus $g=4$ (the dimension of the maximal regular integral manifold, or solution), and the $dx^{i}$ are in involution (so independent variables in the PDE.)  We summarize this well-posed structure by writing 14[0,1,3,5]4+1.  The last integer,~$s_{4}$, is the number of arbitrary functions or gauge freedom in a solution.  There are no Cauchy vectors, so Cartan's test for well-posedness requires that the f\/irst integer equals the sum of the rest.

The Cartan--Poincar\'{e} form for this EDS is
\begin{gather}\label{eq1}
d\Lambda=\theta \wedge \psi.
\end{gather}
That this outer product is indeed closed as we have written follows from a direct computation that $d(\theta  \wedge \psi)=0$.  This in turn is true because
\begin{gather*}
dF \wedge d*F=0
\end{gather*}
which is nothing else than $dF_{ij} \wedge dF_{ij}=0$.
Integrating equation \eqref{eq1} by parts, we take
\begin{gather*}
\Lambda=-*F\wedge dA+\tfrac{1}{2} F \wedge *F.
\end{gather*}
This Cartan form corresponds to the Lagrangian density for Maxwell theory used by f\/ield theo\-rists~\cite[page 98, equation~(9.13)]{Schw}.

The generalization to Yang--Mills theory in 4 dimensions is of course well known.  Its form is nicely exposed when written as an EDS.  Structure constants of an $m$-dimensional Lie group~$\gamma ^{a}_{bc}$, $a,b,\ldots = 1,\dots,m$ are used, taken in a basis where these are completely antisymmetric ($\gamma ^{a}_{bc}=-\gamma ^{a}_{cb}=\gamma ^{b}_{ca}=\gamma ^{c}_{ab}$) and we take $\gamma ^{1}_{23}=1$.  Quadratically nonlinear terms are introduced consistently coupling indexed variables $A^{a}_{i}$  and $F^{a}_{ij}$.  In the following we consider only the case $m=3$, which is $SU(2)$.  We def\/ine 1-forms and 2-forms
\begin{gather*}
A^{a}=A^{a}_{i} dx^{i},\\ 
F^a = \tfrac{1}{2} F^{a}_{ij} dx^{i} \wedge dx^{j},\\
{}*F_a = \tfrac{1}{4} F^{a}_{ij}  dx^{k} \wedge dx^{l} \epsilon^{ij}_{..kl}
\end{gather*}
and now set
\begin{gather*}
\theta^{a} = dA^{a}-F^a+\tfrac{1}{8} \gamma ^{a}_{bc} A^{b} \wedge A^{c},
\\
\psi_{a} = - d*F_{a} - \tfrac{1}{4}\gamma^{c}_{ab } A^{b} \wedge *F_{c}.
\end{gather*}

The EDS in 34 variables is generated by the three 2-forms $\theta^{a}$, three 3-forms $d\theta^{a}$ and three additional 3-forms $\psi_{a}$.  That this EDS is closed follows from calculating
\begin{gather*}
d\psi_{a}=- \tfrac{1}{4} \gamma ^{c}_{ab} A^{b} \wedge \psi_{c} - \tfrac{1}{4} \gamma ^{c}_{ab} \theta^{b} \wedge *F_{c},
\end{gather*}
so $d\psi_{a}$ are in the EDS (the essential identity in this is $\gamma^{c}_{ab } *F_{c} \wedge F^{b} = 0$).

The Cartan--Poincar\'{e} form is now
\begin{gather*}
d \Lambda = \theta^{a} \wedge \psi_{a}
\end{gather*}
and it can be explicitly integrated to take
\begin{gather*}
\Lambda = - *F_{a} \wedge dA^{a} + \tfrac{1}{2} F^{a} \wedge *F_{a}  - \tfrac{1}{8} \gamma^{c}_{ab} A^{a} \wedge  A^{b} \wedge *F_{c}.
\end{gather*}

A possibly new result is a calculation of the Cartan character table.  We f\/ind
\begin{gather*}
34[0,3,9,15]4+3
\end{gather*}
showing this theory to be well-posed, with 3 degrees of gauge freedom.

\section{Generalization to other dimensions}\label{sec3}

All the equations of Section 2 can be taken as they stand in any number of spacetime dimensions, $i, j = 1,\dots,n$.  We need only redef\/ine the dual f\/ield form as an $(n-2)$-form
\begin{gather*}
*F_{a} = \tfrac{1}{2 (n-2)!} F^{a}_{ij} dx^{k} \wedge dx^{l}\wedge\cdots \wedge dx^{n} \epsilon^{ij}_{..kl...n}.
\end{gather*}
The $\psi_{a}$ are now $(n-1)$-forms.  We have explicitly verif\/ied this generalization, and calculated the Cartan characters for theories through $n= 6$, showing the EDS' to be well-posed. These results are summarized in Table~\ref{table1}, and the pattern for higher dimensions is apparent.

\begin{table}[h!]\centering

\caption{Calculated Cartan characters for Maxwell and $SU(2)$-Yang--Mills theories in $n$ dimensions.  The characters $s_{0},\dots,s_{n-1}$ are given in brackets [~~], followed by $n$.  The preceding number is the dimension of the manifold where the EDS is set.  The last number is $s_{n}$, the gauge freedom.}

\vspace{2mm}

\begin{tabular}{l}\label{table1}
Maxwell characters in $n$ dimensions\\[2mm]
\hphantom{2}9[0,2,3]3 + 1\\
14[0,1,3,5]4 + 1\\
20[0,1,2,4,7]5 + 1\\
27[0,1,2,3,5,9]6 + 1\\[4mm]
$SU(2)$-Yang--Mills characters in $n$ dimensions\\[2mm]
21[0,6,9]3 + 3\\
34[0,3,9,15]4 + 3\\
50[0,3,6,12,21]5 + 3\\
69[0,3,6,9,15,27]6 + 3
\end{tabular}
\end{table}

\subsection*{Acknowledgements}

This work 
 was performed while I held a visiting appointment at the Jet Propulsion Laboratory,  California Institute of Technology.

\pdfbookmark[1]{References}{ref}
\LastPageEnding

\end{document}